\documentclass[12pt,reqno,a4paper]{amsart}

\usepackage[mathlines]{lineno}
\AddToHook{env/abstract/begin}{\internallinenumbers}
\AddToHook{env/abstract/end}{\endinternallinenumbers}

\usepackage{geometry}
\usepackage[british]{babel}
\usepackage[utf8]{inputenc}
\usepackage{lmodern}
\usepackage{amsmath,amssymb,amsthm,mathtools}
\usepackage[hidelinks]{hyperref}
\usepackage[capitalise,noabbrev]{cleveref}
\allowdisplaybreaks

\newtheorem{theorem}{Theorem}[section]
\newtheorem{lemma}[theorem]{Lemma}
\newtheorem{proposition}[theorem]{Proposition}
\newtheorem{corollary}[theorem]{Corollary}

\newtheorem{conjecture}[theorem]{Conjecture}
\newtheorem*{claim}{Claim}

\crefname{theorem}{Theorem}{Theorems}
\crefname{lemma}{Lemma}{Lemmas}
\crefname{proposition}{Proposition}{Propositions}
\crefname{corollary}{Corollary}{Corollaries}
\crefname{remark}{Remark}{Remarks}
\crefname{conjecture}{Conjecture}{Conjectures}

\makeatletter
\newskip\nestedblockparskip
\newenvironment{nestedproofblock}
  {%
    \nestedblockparskip=\parskip
    \parskip=0pt%
    \par
    \advance\@totalleftmargin by 1.5em%
    \advance\linewidth by -1em%
    \parshape\@ne\@totalleftmargin\linewidth
    \def\thm@space@setup{%
      \thm@preskip=5pt%
      \thm@postskip=3pt%
    }%
  }
  {%
    \par
  }
\newenvironment{nestedproof}[1][\proofname]
  {%
    \par
    \pushQED{\qed}%
    \normalfont
    \topsep=0pt\relax
    \partopsep=0pt\relax
    \trivlist
    \item[\hskip\labelsep\itshape#1\@addpunct{.}]\ignorespaces
    \parskip=\nestedblockparskip
  }
  {%
    \parskip=0pt%
    \popQED
    \endtrivlist
    \@endpefalse
  }
\makeatother

\makeatletter
\let\origsection=\section 
\def\section{\@ifstar{\origsection*}{\mysection}}
\def\mysection{\@startsection{section}{1}\z@{.7\linespacing\@plus\linespacing}{.5\linespacing}{\normalfont\scshape\centering\S}}
\makeatother

\AddToHook{env/lemma/begin}{%
  \crefalias{theorem}{lemma}%
}

\AddToHook{env/proposition/begin}{%
  \crefalias{theorem}{proposition}%
}

\AddToHook{env/corollary/begin}{%
  \crefalias{theorem}{corollary}%
}
\AddToHook{env/conjecture/begin}{%
  \crefalias{theorem}{conjecture}%
}
\AddToHook{env/remark/begin}{%
  \crefalias{theorem}{remark}%
}

\newcommand{\lir}{\chi'_{\mathrm{irr}}}

\hypersetup{
  pdftitle={An 18-colour bound for locally irregular decompositions},
  pdfauthor={Carla Negri Lintzmayer, Guilherme Oliveira Mota, Maycon
  Sambinelli, and Vinicius Fernandes dos Santos},
  pdfsubject={Locally irregular edge-decompositions},
  pdfkeywords={locally irregular graph,
    edge-decomposition, locally irregular chromatic index}
}

\title{An 18-colour bound for locally irregular decompositions}

\author[C.~N.~Lintzmayer]{Carla Negri Lintzmayer}
\address{Centro de Matem\'atica, Computa\c{c}\~ao e Cogni\c{c}\~ao,
Universidade Federal do ABC, Santo Andr\'e, Brazil}
\email{carla.negri@ufabc.edu.br}

\author[G.~O.~Mota]{Guilherme Oliveira Mota}
\address{Departamento de Ci\^encia da Computa\c{c}\~ao,\newline Instituto de
Matem\'atica, Estat\'{\i}stica e Ci\^encia da Computa\c{c}\~ao, Universidade de
S\~ao Paulo, S\~ao Paulo, Brazil}
\email{mota@ime.usp.br}

\author[M.~Sambinelli]{Maycon Sambinelli}
\address{Centro de Matem\'atica, Computa\c{c}\~ao e Cogni\c{c}\~ao,
Universidade Federal do ABC, Santo Andr\'e, Brazil}
\email{m.sambinelli@ufabc.edu.br}

\author[V.~F.~dos Santos]{Vinicius Fernandes dos Santos}
\address{Departamento de Ci\^encia da Computa\c{c}\~ao, Universidade Federal de
Minas Gerais, Belo Horizonte, Brazil}
\email{viniciussantos@dcc.ufmg.br}

\subjclass[2020]{Primary 05C15; Secondary 05C70}
\keywords{locally irregular graph, edge-decomposition, locally irregular
chromatic index}
\date{}

\begin{document}

\begin{abstract}
  A graph is \emph{locally irregular} if adjacent vertices have distinct
  degrees. A graph~$G$ is \emph{decomposable} if its edge set can be decomposed
  into locally irregular graphs, and its \emph{locally irregular chromatic
  index} $\lir(G)$ is the least number of graphs in such a decomposition. We
  prove that $\lir(G)\le18$ for every decomposable graph~$G$, improving the
  previous bound of~$220$.
\end{abstract}

\maketitle

\section{Introduction}

Throughout this paper, all graphs are finite and simple. A graph is
\emph{locally irregular} if the endvertices of every edge have distinct
degrees. 
An edge-colouring is locally irregular if every colour class forms a locally
irregular graph. 
A graph~$G$ is \emph{decomposable} if it admits such a colouring, and $\lir(G)$
is the least number of colours in such a colouring. 
The family of non-decomposable graphs was characterised by Baudon, Bensmail,
Przyby{\l}o, and Wo{\'z}niak~\cite{BaBePrWo15}, who conjectured that every
decomposable graph~$G$ has $\lir(G)\le3$.
However, Sedlar and \v{S}krekovski~\cite{SeSk21} exhibited the \emph{bow-tie},
obtained from an edge by attaching two triangles at each end, with all new
vertices distinct, and proved that its locally irregular chromatic index
is~$4$.
They later proposed the following conjecture~\cite{SeSk24}.

\begin{conjecture}\label{conj:main}
  If $G$ is a connected decomposable graph that is not the bow-tie, then
  $\lir(G) \le 3$.
\end{conjecture}

Baudon, Bensmail, and Sopena~\cite{BaBeSo15} proved that deciding whether a
given graph admits a locally irregular edge-colouring with at most two colours
is NP-complete, even for graphs of maximum degree at most~$6$. 
In contrast, for trees they gave a linear-time algorithm that detects the
non-decomposable case and otherwise computes~$\lir(T)$. 
Przyby{\l}o~\cite{Pr16} proved that every decomposable graph of sufficiently
large minimum degree satisfies $\lir(G) \le 3$. 
Bensmail, Merker, and Thomassen~\cite{BeMeTh17} proved that $\lir(G) \le 328$
for any decomposable graph~$G$, which was later reduced to~$220$ by Lu\v{z}ar,
Przyby{\l}o, and Sot\'ak~\cite{LuPrSo18}. 
For results on particular classes of graphs, see, for
example,~\cite{BeDrNi20,LeLiShZh22,LiMoSa21,LuMaRiSoSrSt23,Pr24,SeSk24}. 
Our main result improves the general upper bound from~$220$ to~$18$.

\begin{theorem}\label{thm:global18}
  Every decomposable graph~$G$ satisfies $\lir(G) \le 18$.
\end{theorem}

A combination of Theorem~2.3 and Lemma~3.2 of Bensmail, Merker, and
Thomassen~\cite{BeMeTh17} with Corollary~10 of Lu\v{z}ar, Przyby{\l}o, and
Sot\'ak~\cite{LuPrSo18} gives the following result, which may be of independent
interest.
For completeness, we show the proof in Section~\ref{sub:sec:chromatic_bound}.

\begin{theorem}\label{thm:chromatic-bound}
  Let~$G$ be a decomposable connected graph, and let $\tau \in \{0,1\}$ satisfy
  $\tau \equiv |E(G)| \pmod 2$. 
  Then $\lir(G) \le 4\chi(G) - 2 + \tau$.
\end{theorem}

Given a positive integer~$k$ and a graph $G$, a \emph{$k$-labelling} of~$V(G)$
is simply a map $\varphi\colon V(G)\to\{1,\ldots,k\}$. 
An edge is \emph{crossing} when its endpoints have distinct labels and
\emph{monochromatic} otherwise. 
The spanning subgraphs of~$G$ consisting of these two types of edges are the
\emph{crossing graph} and the \emph{monochromatic graph}, respectively. 
The labelling is \emph{crossing-optimal} if it maximises the number of crossing
edges.

\subsection*{Overview of the proof}
Since colours may be reused in distinct components, it is enough to consider a
connected graph~$G$. 
Fix a crossing-optimal $6$-labelling of~$V(G)$, and let~$C$ and~$F$ be its
crossing and monochromatic graphs. 
A switching argument shows that~$C$ is connected (\Cref{lem:switching}). 
If $F=\emptyset$, our $16$-colour bound for decomposable $6$-colourable graphs
applies (\Cref{lem:even-fifteen}).

Assume $F\neq\emptyset$ and group the six labels into three pairs
$B_1, B_2, B_3$. 
For $i \in \{1,2,3\}$, let~$J_i$ be the subgraph of~$F$ induced by the vertices
with a label in~$B_i$; these three graphs edge-partition~$F$.
For each~$J_i$ with an edge, orient it so that a subset of crossing edges can
later be added to make the union locally irregular, with the necessary 
edge-count parity (\Cref{lem:general-core}). 
Four crossing edges can then be assigned to each oriented monochromatic edge,
one to each label outside its pair, with all assigned edges distinct and such
that deleting them leaves~$C$ connected (\Cref{cor:support}). 
For each~$J_i$, choose a subset~$Q_i$ of its assigned edges so that $J_i \cup
Q_i$ is locally irregular. 
Since $F \neq \emptyset$, choose~$i_0$ with $E(J_{i_0}) \neq \emptyset$ and
whenever~$i \neq i_0$ with $E(J_i) \neq \emptyset$, assure that~$e(Q_i)$ is even,
and that $e(Q_{i_0}) \equiv e(C) \pmod 2$. 
The unused crossing graph then has an even number of edges, it remains
connected, and it is $6$-colourable, so it needs at most~$15$ further classes
(\Cref{lem:even-fifteen}). 
The total number of colours is thus at most $3+15=18$.

\section{Structural information by orientations}

Auxiliary multigraphs are finite and loopless, but may have parallel edges.
For a graph or auxiliary multigraph~$H$, put $e(H) := |E(H)|$ and write~$d_H(v)$
for the degree of~$v$. 
If~$\vec H$ is an orientation of~$H$, write $d^-_{\vec H}(v)$ for its indegree.
For a set~$Y$ of vertices, write~$d_H(v,Y)$ for the number of edges from~$v$
to~$Y$, and let~$H[Y]$ denote the subgraph induced by~$Y$.
For an edge set~$A$, let $H-A$ be the spanning subgraph with edge set
$E(H) \setminus A$. 
We use the following elementary fact. 
Its even-orientation special case is standard; see, for example, Cannon,
Ishaque, and T\'oth~\cite{CaIsTo12}. 
We include a short proof for completeness.

\begin{samepage}
\begin{lemma}\label{lem:parity-orientation}
  Let~$L$ be a connected loopless multigraph. 
  Then
  \begin{enumerate}
    \renewcommand{\theenumi}{\roman{enumi}}
    \item For every map $\eta \colon V(L) \to \{0,1\}$ and every spanning tree~$T$
    of~$L$, every orientation of $L-E(T)$ extends to an orientation $\vec L$
    of~$L$ satisfying $d^-_{\vec L}(v) \equiv \eta(v) \pmod 2$ for all~$v$ if and
    only if $\sum_{v \in V(L)} \eta(v) \equiv e(L) \pmod 2$.
    \item There is an orientation of $L$ in which all indegrees are even,
    except possibly one, which may have indegree~$1$.
  \end{enumerate}
\end{lemma}
\end{samepage}
\begin{proof}
  Necessity in~(i) follows from $\sum_{v \in V(L)} d^-_{\vec L}(v) = e(L)$.
  Conversely, suppose that $\sum_{v \in V(L)} \eta(v) \equiv e(L) \pmod 2$, fix
  an arbitrary orientation of~$L-E(T)$, root~$T$ at~$r$, and process the
  non-root vertices with children preceding their parents. 
  At each non-root vertex, its parent edge is the only unoriented incident
  edge; orient it to obtain the prescribed parity. 
  The resulting orientation~$\vec L$ also has the prescribed parity at~$r$,
  since
  \[
    d^-_{\vec L}(r) \equiv e(L) - \sum_{v \neq r} \eta(v) \equiv \eta(r) \pmod 2 \,.
  \]

  For~(ii), the assertion is immediate if $e(L) = 0$. 
  Otherwise, choose a leaf~$\ell$ of a spanning tree~$T$, orient all non-tree
  edges at~$\ell$ away from~$\ell$, and orient the other non-tree edges
  arbitrarily.
  Apply~(i) with $\eta(\ell) \equiv e(L) \pmod 2$ and $\eta(v) = 0$ elsewhere,
  and let~$\vec L$ be the resulting orientation. 
  Since only the tree edge at~$\ell$ can enter~$\ell$, we have $d^-_{\vec
  L}(\ell) \in \{0,1\}$, while every other vertex has even indegree.
\end{proof}

The following result is due to Akbari, Dalirrooyfard, Ehsani, Ozeki, and
Sherkati~\cite[Lemma~11 in the authors' preprint]{AkDaEhOzSh20}.
Reversing all orientations gives the indegree version below.
We include a short proof to make explicit that it also holds for multigraphs.

\begin{lemma}\label{lem:lists} 
  Let~$N$ be a bipartite loopless multigraph with finite forbidden sets of
  non-negative integers~$F(v)$ for each vertex~$v$.
  If~$N$ has an orientation~$D$ satisfying $d^-_D(v) \ge |F(v)|$ at every
  vertex~$v$, then~$N$ has an orientation~$O$ satisfying $d^-_O(v) \notin F(v)$
  at every~$v$.
\end{lemma}
\begin{proof}
  Fix an orientation~$D$ of~$N$ such that $d^-_D(x) \ge |F(x)|$ for every~$x
  \in V(N)$, and proceed by induction on~$e(N)$. 
  If no~$F(x)$ contains either~$0$ or~$d_N(x)$, orient every edge from one
  fixed bipartition class to the other. 
  Every indegree is then either~$0$ or~$d_N(x)$, so this orientation
  avoids~$F$.

  Otherwise, choose~$u$ such that $0 \in F(u)$ or $d_N(u) \in F(u)$. 
  Since $d^-_D(u) \ge |F(u)| \ge 1$, some edge~$vu$ is directed from~$v$ to~$u$
  in~$D$.
  Put $N' := N - vu$ and let~$D'$ be the restriction of~$D$ to~$N'$.

  First suppose that $0 \in F(u)$. 
  Define
  \[
    F'(u) := \{j-1 \colon j \in F(u) \setminus \{0\}\},
            \qquad F'(x) := F(x) \quad (x \neq u) \, .
  \]
  We have $d^-_{D'}(u) = d^-_D(u) - 1 \ge |F'(u)|$, while $d^-_{D'}(x) =
  d^-_D(x) \ge |F'(x)|$ for every~$x \neq u$.
  By induction, $N'$ has an orientation~$O'$ satisfying $d^-_{O'}(x) \notin
  F'(x)$ for every~$x$.
  Restore~$vu$ directed from~$v$ to~$u$, obtaining an orientation~$O$ of~$N$.
  Then $d^-_O(u) = d^-_{O'}(u) + 1 > 0$.
  If $d^-_O(u) = j \in F(u) \setminus \{0\}$, then $d^-_{O'}(u) = j-1 \in
  F'(u)$, a contradiction. 
  For every~$x \neq u$, the restored edge is not directed into~$x$, so
  $d^-_O(x) = d^-_{O'}(x) \notin F'(x) = F(x)$.
  Thus~$O$ avoids~$F$.

  We may therefore assume that $0 \notin F(u)$ and $d_N(u) \in F(u)$.
  Define
  \[
    F'(u) := F(u) \setminus \{d_N(u)\}, \qquad
    F'(v) := \{j-1 \colon j \in F(v) \setminus \{0\}\} \, ,
  \]
  and let $F'(x) := F(x)$ for $x \notin \{u,v\}$.
  Again, $d^-_{D'}(u) = d^-_D(u)-1 \ge |F'(u)|$; moreover, $d^-_{D'}(v) =
  d^-_D(v) \ge |F(v)| \ge |F'(v)|$, while $d^-_{D'}(x) = d^-_D(x) \ge |F(x)| =
  |F'(x)|$ for every $x \notin \{u,v\}$.
  By induction, $N'$ has an orientation~$O'$ satisfying $d^-_{O'}(x) \notin
  F'(x)$ for every~$x$.
  Restore~$vu$ directed from~$u$ to~$v$, obtaining an orientation~$O$ of~$N$.
  Now $d^-_O(u) = d^-_{O'}(u) \notin F'(u)$ and $d^-_O(u) \le d_{N'}(u) =
  d_N(u)-1$, so $d^-_O(u) \notin F(u)$.
  Moreover, $d^-_O(v) = d^-_{O'}(v) + 1 > 0$; if $d^-_O(v) = j \in F(v)
  \setminus \{0\}$, then $d^-_{O'}(v) = j-1 \in F'(v)$, a contradiction.
  All other indegrees are unchanged, and hence~$O$ avoids~$F$.
\end{proof}

By Hall's theorem, the hypothesis of \Cref{lem:lists} holds if and only if
\begin{equation}\label{eq:hall}
  \sum_{v \in S} |F(v)| \le e_N(S) + |\partial_N(S)|
    \qquad\text{for every }S \subseteq V(N) \, ,
\end{equation}
where $e_N(S)$ counts edges with both ends in~$S$, and $\partial_N(S)$ is the
set of edges with exactly one end in~$S$.
Indeed, match $|F(v)|$ copies of each vertex~$v$ to distinct incident edges,
and direct those edges towards their assigned vertices.

\begin{lemma}\label{lem:marked}
  Let~$M$ be a bipartite loopless multigraph, and prescribe a common
  bipartition direction on a matching~$P$ of at most two edges.
  Let~$k$ be even, with $k \ge 2$ if $|P| \le 1$ and $k \ge 4$ otherwise.
  The prescription extends to an orientation whose indegrees are even or~$1$,
  and in which no head of~$P$ has indegree~$k$.
\end{lemma}
\begin{proof}
  Call the heads of~$P$ marked and all other vertices ordinary.
  Put $N := M-P$.
  Once the prescribed arcs of~$P$ are restored, the total indegree of a marked
  head is its indegree in~$N$ plus~$1$, whereas the total indegree of an
  ordinary vertex is its indegree in~$N$.
  Accordingly, for an ordinary vertex~$v$, define
  \[
    F(v) := \{r \in \{3,4,\ldots,d_N(v)\} \colon r \text{ is odd }\} \, ,
  \]
  and, for a marked vertex~$v$, define
  \[
    F(v) := \{r \in \{0,1,\ldots,d_N(v)\} \colon r \text{ is positive and even,
    or } r = k-1\}.
  \]
  Put $\delta(v) := 2|F(v)| - d_N(v)$.
  Since the prescribed edges have a common bipartition direction, all marked
  heads lie in the same bipartition class.

  At an ordinary nonisolated vertex, $\delta(v)$ is~$-1$ when $d_N(v)$ is odd
  and~$-2$ when~$d_N(v)$ is even.
  At a marked head with $d_N(v) \ge k-1$, the value $\delta(v)$ is~$1$ when
  $d_N(v)$ is odd and~$2$ when $d_N(v)$ is even; at any other marked head it
  is~$-1$ or~$0$.
  An isolated ordinary vertex has $\delta(v) = 0$.

  If~\eqref{eq:hall} holds, the preceding Hall criterion and \Cref{lem:lists}
  give an orientation of~$N$ avoiding every~$F(v)$.
  After the prescribed arcs of~$P$ are restored, every total indegree is even
  or~$1$, and no marked head has total indegree~$k$.

  We may therefore assume that~\eqref{eq:hall} fails.
  The difference $\sum_{v \in S} |F(v)| - e_N(S) - |\partial_N(S)|$ is additive
  over the components of~$N[S]$, so there is a set~$S$ such that~$N[S]$ is
  connected and this difference is positive. 
  Since $\sum_{v \in S} d_N(v) = 2e_N(S) + |\partial_N(S)|$, this is equivalent
  to
  \begin{equation}\label{eq:deficiency}
    \sum_{v \in S} \delta(v) - |\partial_N(S)|
      = 2 \left(\sum_{v \in S} |F(v)| - e_N(S) - |\partial_N(S)| \right) \ge 2 \, .
  \end{equation}
  A singleton cannot satisfy~\eqref{eq:deficiency}: in that case
  $|\partial_N(S)| = d_N(v)$, while the values above and the assumptions on~$k$
  give $\delta(v) \le d_N(v)$. 
  If a connected set with at least two vertices contains at most one marked
  head~$v$ satisfying $d_N(v) \ge k-1$, its $\delta$-sum is at most~$1$: that
  head contributes at most~$2$ and, if present, has an ordinary neighbour
  contributing at most~$-1$, while all other vertices contribute at most~$0$.
  Thus~$S$ contains both marked heads, say~$b_1,b_2$, each satisfying
  $d_N(b_i) \ge k-1$, and $k \ge 4$.

  Let~$q$ be the number of ordinary vertices in~$S$.
  Since the two heads lie in the same bipartition class and $N[S]$ is
  connected, $q \ge 1$.
  Each head contributes at most~$2$ to the $\delta$-sum, and each ordinary
  vertex contributes at most~$-1$. 
  Hence~\eqref{eq:deficiency} gives $2 \le 4 - q - |\partial_N(S)|$.
  Thus either $q = 1$ and $|\partial_N(S)| \le 1$, or $q = 2$ and
  $|\partial_N(S)| = 0$.
  In the latter case, $S$ is the vertex set of a component of~$N$, both heads
  have even degree in~$N$, and both ordinary vertices have odd degree in~$N$.

  If $q = 1$, let~$a$ be the unique ordinary vertex in~$S$.
  Direct every edge of~$N[S]$ towards~$a$ and direct the possible boundary edge
  into~$S$.
  If the resulting indegree of~$a$ in~$N$ is odd, reverse one edge of~$N[S]$
  towards a head not incident with the boundary.
  Such an edge exists because~$N[S]$ is connected. 
  The indegree of~$a$ is then even and, after the prescribed arcs of~$P$ are
  restored, each marked head has total indegree~$1$ or~$2$.

  If $q = 2$, the underlying simple graph of~$N[S]$ has a perfect matching.
  Direct every edge of~$N[S]$ towards the ordinary vertices and reverse the
  direction of one edge corresponding to each edge of this matching. 
  Each ordinary vertex then has even indegree, since its degree in~$N$ is odd,
  and each marked head has one incoming edge from~$N[S]$ and one from~$P$,
  hence total indegree~$2$.

  There are no marked heads outside~$S$. 
  The possible boundary edge was directed into~$S$, so it contributes no
  indegree outside~$S$.
  Apply part~(ii) of \Cref{lem:parity-orientation} separately to the nontrivial
  components of $N[V(N) \setminus S]$, orienting their internal edges so that
  every indegree is even or~$1$.
  Together, these orientations and those already fixed inside~$S$, on the
  possible boundary edge, and on~$P$ give every vertex of~$M$ indegree even
  or~$1$ and give neither marked head indegree~$k$.
\end{proof}

The next lemma shows how to combine a monochromatic subgraph with suitable
crossing edges to obtain a locally irregular graph, while prescribing the
parity of the number of crossing edges used.

\begin{lemma}\label{lem:general-core}
  Let~$J$ be a graph with an edge and vertex set~$X$.
  There exists an orientation~$\vec J$ of~$J$ such that every bipartite graph
  $B=(X,Y_1\mathbin{\dot\cup}Y_2; E)$ satisfying
  \[
    d_B(v,Y_i) = 2d^-_{\vec J}(v) 
        \qquad \text{for every } v \in X,\ i \in \{1,2\}
  \]
  contains, for each $s \in \{0,1\}$, a subgraph~$Q_s$ such that $e(Q_s) \equiv
  s \pmod 2$ and $J \cup Q_s$ is locally irregular.
\end{lemma}
\begin{proof}
  Choose a vertex~$x$ of maximum degree~$\Delta$ in~$J$. 
  Among the orientations of~$J$ in which every edge incident with~$x$ is
  directed towards~$x$, choose~$\vec J$ to maximise $\sum_{v \in X} b(v)^2$,
  where
  \[
    b(v) := d_J(v) + 2d^-_{\vec J}(v) \, .
  \]
  Reversing an arc $(u,v)$ not incident with~$x$ changes this sum by
  $4(b(u)-b(v)+2)$.
  Hence maximality gives $b(v) \ge b(u) + 2$.
  The same holds on every arc $(u,x)$, since $b(u) \le 3d_J(u) - 2 \le 3\Delta
  - 2 = b(x)-2$.

  Let~$B$ satisfy the displayed degree condition. 
  For any integers $0 \le q(v) \le 2d^-_{\vec J}(v)$, choose~$q(v)$ edges from
  each~$v$ to each of~$Y_1,Y_2$, and pair them across the two sets. 
  Replace each pair $vy,vy'$ by an auxiliary edge~$yy'$.
  The resulting multigraph~$M$ is bipartite and loopless. 
  Any orientation $\vec M$ of~$M$ selects a subgraph~$Q\subseteq B$ by
  retaining, from each pair, the edge incident with its auxiliary head. 
  Thus
  \[
    d_Q(v) = q(v) \quad (v \in X), \qquad
    d_Q(y) = d^-_{\vec M}(y) \quad (y \in Y_1 \cup Y_2) \, .
  \]

  First put $q(v) := 2d^-_{\vec J}(v) - \varepsilon(v)$, where $\varepsilon(v)
  = 1$ if $d_J(v)$ is even and $d^-_{\vec J}(v) > 0$, and $\varepsilon(v) = 0$
  otherwise.
  Apply part~(ii) of \Cref{lem:parity-orientation} separately to the nontrivial
  components of~$M$, and let~$\vec M$ be the union of the resulting
  orientations.
  Let~$Q$ be the subgraph selected by~$\vec M$.
  In $H := J \cup Q$, every~$v \in X$ has degree $b(v) - \varepsilon(v)$.
  Hence, for every arc $(u,v)$ of~$\vec J$,
  \[
    d_H(v) - d_H(u)
      =b(v) - b(u) - \varepsilon(v) + \varepsilon(u) \ge 1 \, .
  \]
  Moreover, every selected edge joins a vertex~$v \in X$ with $d^-_{\vec J}(v)
  > 0$, for which $d_H(v)$ is odd and at least~$3$, to a vertex outside~$X$
  whose degree is even or~$1$.
  Thus~$H$ is locally irregular.

  For the other parity, change only~$q(x)$: use~$1$ if~$\Delta$ is odd
  and~$2$ if~$\Delta$ is even.
  The auxiliary edges arising at~$x$ form a matching~$P$, since~$B$ is simple.
  Prescribe them towards~$Y_2$ and apply \Cref{lem:marked} with $k := \Delta +
  q(x)$.
  This is even and exceeds~$\Delta$; moreover, $k \ge 4$ when $|P| = 2$.
  Let~$Q$ be the subgraph selected by the resulting orientation.
  In the resulting graph~$H := J \cup Q$, we have $d_H(x) = k$.
  For every arc $(u,v)$ of~$\vec J$ not incident with~$x$,
  \[
    d_H(v) - d_H(u)
      = b(v) - b(u) - \varepsilon(v) + \varepsilon(u) \ge 1 \, .
  \]
  On an edge
  $xu \in E(J)$, if $d^-_{\vec J}(v) > 0$, then~$d_H(u)$ is odd, whereas~$k$ is
  even; if $d^-_{\vec J}(v) = 0$, then $d_H(u) = d_J(u) \le \Delta < k$.
  Every selected edge at~$x$ ends at a head of~$P$, whose degree is different
  from~$k$ by \Cref{lem:marked}. 
  Finally, every selected edge not incident with~$x$ joins an odd-degree vertex
  of~$X$ with degree at least~$3$ to a vertex outside~$X$ of even degree or
  degree~$1$.
  Thus~$H$ is locally irregular.

  The original value of~$q(x)$ was~$2\Delta$ for odd~$\Delta$ and $2\Delta-1$
  for even~$\Delta$.
  The change therefore reverses the parity of $e(Q) = \sum_{v \in X} q(v)$.
  Choosing the construction with edge-count parity~$s$ gives the subgraph~$Q_s$
  for each $s \in \{0,1\}$.
\end{proof}

\section{Chromatic bounds}

\subsection{A general chromatic bound}
\label{sub:sec:chromatic_bound}

A bipartite graph is \emph{semi-even} if, in every nontrivial component, one
bipartition class consists entirely of even-degree vertices; this class may be
chosen independently in different components. 
Such a graph is decomposable since it has an even total number of edges. 

\begin{proposition}\label{prop:semi-even}
  Every semi-even forest has locally irregular chromatic index at most~$2$
  \cite[Lemma~3.2]{BeMeTh17}.
  Every semi-even bipartite graph has locally irregular chromatic index at
  most~$4$ \cite[Corollary~10]{LuPrSo18}.
\end{proposition}

We also use the following result of Bensmail, Merker, and
Thomassen~\cite[Theorem~2.3]{BeMeTh17}.

\begin{lemma}\label{lem:odd-reduction}
  Let~$G$ be a connected decomposable graph.
  If $|E(G)|$ is odd, then~$G$ contains a locally irregular subgraph~$L$ such
  that each component of $G-E(L)$ contains an even number of edges.
\end{lemma}

The next argument extends to multipartite graphs the semi-even decomposition
method used for bipartite graphs~\cite{BeMeTh17,LuPrSo18}.
For a positive integer~$k$, write $[k]:=\{1,\ldots,k\}$.

\begin{lemma}\label{lem:partite-residual}
  Let~$H$ be a graph in which every component has an even number of edges.
  Then~$H$ can be edge-decomposed into one possibly empty semi-even forest and
  at most $\chi(H)-1$ semi-even bipartite graphs.
  Consequently,
  \[
    \lir(H)\le4\chi(H)-2.
  \]
\end{lemma}
\begin{proof}
  Put $k:=\chi(H)$ and fix a partition $V_1 \cup \cdots \cup V_k$ of~$V(H)$
  into independent sets. In every nontrivial component of~$H$, choose a
  spanning tree, and let~$T$ be the union of these trees.
  Orient every edge of $H-E(T)$ incident with~$V_1$ away from~$V_1$, and orient
  the other edges of $H-E(T)$ arbitrarily.
  Since every component has an even number of edges, part~(i) of
  \Cref{lem:parity-orientation}, applied separately in each nontrivial
  component with even indegree prescribed at every vertex, extends the
  orientation to all edges of~$T$.

  For $i \in [k]$, let~$M_i$ be the spanning subgraph containing the edges
  oriented towards~$V_i$. These graphs edge-partition~$H$. Each~$M_i$ is
  bipartite with even degrees in~$V_i$, and hence semi-even.
  Moreover, $M_1\subseteq T$ by the choice of the non-tree orientations.
  Thus \Cref{prop:semi-even} gives
  $\lir(H)\le2+4(k-1)=4\chi(H)-2$.
\end{proof}

\begin{proof}[Proof of \Cref{thm:chromatic-bound}]
  If $\tau=1$, \Cref{lem:odd-reduction} gives a locally irregular subgraph~$L$
  such that every component of $G-E(L)$ has an even number of edges; if
  $\tau=0$, let~$L$ be empty. Put $G':=G-E(L)$. By
  \Cref{lem:partite-residual},
  $\lir(G)\le\tau+\lir(G') \le\tau+4\chi(G')-2 \le4\chi(G)-2+\tau$.
\end{proof}

Applying \Cref{thm:chromatic-bound} separately to the components and reusing
the same set of colours gives the best-known bounds of~$7$ for decomposable
bipartite graphs and of~$15$ for decomposable planar
graphs~\cite{LuPrSo18,BeDrNi20}.
If every component has an even number of edges, \Cref{lem:partite-residual}
also matches the known bounds of~$6$ and~$14$, respectively.

\subsection{The six-colourable case}

\begin{lemma}\label{lem:even-fifteen}
  Let~$G$ be a graph with $\chi(G) \le 6$.
  If every component of~$G$ has an even number of edges, then $\lir(G) \le 15$.
  Without the parity assumption, every decomposable graph~$G$ with $\chi(G) \le
  6$ satisfies $\lir(G) \le 16$.
\end{lemma}
\begin{proof}
  The result is trivial if~$G$ is edgeless.
  For the first assertion, it suffices to prove the bound separately for each
  nontrivial component and reuse the same colours. 
  Thus we may assume that~$G$ is connected and has at least one edge.

  Identify the six available labels with~$\mathbb Z_6$.
  Among all proper maps $\psi \colon V(G) \to \mathbb Z_6$, choose one
  maximising the number of edges whose endpoint labels differ by~$\pm 1$.
  Define the spanning subgraphs $U,T,O$ (unit, two-step, and opposite) by
  \begin{align*}
    E(U) & := \{xy \in E(G) \colon |\psi(x) - \psi(y)| \in \{1,5\}\},\\
    E(T) & := \{xy \in E(G) \colon |\psi(x) - \psi(y)| \in \{2,4\}\},\\
    E(O) & := \{xy \in E(G) \colon |\psi(x) - \psi(y)| = 3\} \, .
  \end{align*}
  These three graphs edge-partition~$G$. 
  Moreover, $O \cup U$ is bipartite under the labels modulo~$2$, while $T \cup
  U$ is $3$-colourable under the labels modulo~$3$.
  Let~$A$ and~$B$ be the sets of vertices with even and odd labels,
  respectively.
  Every component of~$T$ lies in one of these two sets.

  \begin{nestedproofblock}
  \begin{claim}
    The graph~$U$ is connected.
  \end{claim}
  \begin{nestedproof}[Proof of the claim]
    Suppose otherwise, and let~$W$ be the vertex set of a component of~$U$.
    Since~$G$ is connected, at least one edge joins~$W$ to its complement.
    Write each such edge as~$xy$ with $x \in W$, and let
    \[ 
      D_0 := \{|\psi(x) - \psi(y)| \colon xy \in E(G),\ x \in W,\ y \notin W\} \,.
    \]
    None of these edges belongs to~$U$, so $\emptyset \neq D_0 \subseteq
    \{2,3,4\}$.
    Choose $r = 1$ if $4 \in D_0$, $r = -1$ if $4 \notin D_0$ and $2 \in D_0$,
    and $r = 2$ otherwise, and let $D_0^r = \{\ell + r \colon \ell \in D_0\}$.
    In each case, $0 \notin D_0^r$ and $D_0^r \cap \{1,5\} \neq \emptyset$.
    Thus the new map remains proper, preserves every old $U$-edge, and gains a
    new one, contradicting the choice of~$\psi$.
  \end{nestedproof}
  \end{nestedproofblock}

  Put $D := \{a \in A \colon d_O(a) \text{ is odd}\}$.
  Call a component~$K$ of~$T[A]$ \emph{bad} if $e(K)$ is odd and $V(K)
  \subseteq D$.
  In each bad~$K$, choose an edge~$f_K$ such that every component of $K-f_K$
  has an even number of edges: take a non-bridge if~$K$ contains a cycle, and a
  leaf edge otherwise.
  Choose an end~$v_K$ of~$f_K$ and an edge~$v_Kb_K \in E(O)$, which exists
  because $v_K \in D$.

  Group the chosen $O$-edges by their endpoint in~$B$.
  Put all of them in a graph~$P$ and, for each group of odd size, also put
  in~$P$ the edge~$f_K$ associated with one of the chosen $O$-edges in that
  group.
  The subgraphs of~$P$ corresponding to distinct groups are vertex-disjoint:
  vertices contributed in~$A$ by different choices of~$K$ lie in different
  components of~$T[A]$, and different groups have different centres in~$B$.
  An even group is an even star.
  An odd group is an odd star with one extra edge attached at a leaf; this is a
  two-edge path when the group has size~$1$. 
  For a larger odd group, the centre has odd degree at least~$3$, the extended
  leaf has degree~$2$, and all other vertices have degree~$1$. 
  Thus~$P$ is locally irregular and every component of~$P$ has an even number
  of edges.

  Put $O' := O - E(P)$, $T' := T - E(P)$ and $D' := \{a \in A \colon d_{O'}(a)
  \text{ is odd}\}$.
  Each~$v_K$ loses exactly one incident $O$-edge, so $d_{O'}(v_K)$ is even and
  $v_K \notin D'$.
  If $f_K \notin E(P)$, then~$K$ remains a component of~$T'[A]$ with an odd
  number of edges and contains $v_K \notin D'$.
  If $f_K \in E(P)$, every component of $K-f_K$ has an even number of edges.
  Every other component of~$T[A]$ with an odd number of edges was not bad and
  hence contained a vertex outside~$D$; neither its edges nor the $O$-degrees
  of its vertices were changed. 
  Consequently, every component of~$T'[A]$ with an odd number of edges contains
  a vertex outside~$D'$. 
  Since both $e(G)$ and $e(P)$ are even, so is
  \[
    e(U \cup T' \cup O') = e(G) - e(P) \,.
  \]

  In each component of~$T'$ with an odd number of edges, choose a vertex
  outside~$D'$. 
  Such a vertex exists by the previous paragraph for components contained
  in~$A$, while every vertex of a component contained in~$B$ lies outside~$D'$.
  Apply part~(i) of \Cref{lem:parity-orientation} separately to the components
  of~$T'$, prescribing odd indegree at the chosen vertex of each component with
  an odd number of edges and even indegree at every other vertex. 
  Denote the resulting orientation by~$\vec T'$. 
  In particular, $d^-_{\vec T'}(a)$ is even for every $a \in D'$.

  Choose a spanning tree~$W$ of~$U$ and orient the edges of $U-E(W)$
  arbitrarily. 
  Define $\eta \colon V(G) \to \{0,1\}$ by requiring
  \[
    \eta(v) \equiv
    \begin{cases}
      1-d^-_{\vec T'}(v), & v\in D',\\
      -d^-_{\vec T'}(v), & v\notin D'
    \end{cases}
    \pmod 2 \,.
  \]
  Since $O'$ is bipartite with parts~$A$ and~$B$,
  \[
    |D'| \equiv \sum_{a \in A} d_{O'}(a) = e(O') \pmod 2 \,.
  \]
  Together with the evenness of $e(U \cup T' \cup O')$, this gives
  \[
    \sum_{v \in V(G)} \eta(v)
      \equiv |D'| - \sum_{v \in V(G)} d^-_{\vec T'}(v)
      = |D'| - e(T')
      \equiv e(U) \pmod 2 \, .
  \]
  Part~(i) of \Cref{lem:parity-orientation} therefore extends the chosen
  orientation to an orientation~$\vec U$ of~$U$ satisfying $d^-_{\vec U}(v)
  \equiv \eta(v) \pmod 2$ for every vertex~$v$. 
  Hence
  \[
    d^-_{\vec U}(v) + d^-_{\vec T'}(v)
  \]
  is odd precisely when $v \in D'$.
  Moreover, if $a \in D'$, then $d^-_{\vec T'}(a)$ is even, so $d^-_{\vec U}(a)$
  is odd and therefore positive.

  Choose one edge of~$U$ directed into each vertex of~$D'$, and let~$S$
  be their set.
  Every edge of~$U$ has exactly one endpoint in~$A$, so the chosen edges are
  distinct and
  \[
    d_S(a) =
    \begin{cases}
      1, & a \in D',\\
      0, & a \in A \setminus D' \, .
    \end{cases}
  \]
  Thus $O'\cup S$ is bipartite and has even degree at every vertex of~$A$, so
  \Cref{prop:semi-even} gives at most four colours.

  Let~$\vec R$ be the orientation of $R := (U - S) \cup T'$ inherited from
  $\vec U$ and~$\vec T'$. 
  At each vertex of~$D'$, deleting~$S$ removes exactly one incoming edge, while
  no edge of~$S$ enters a vertex outside $D'$.
  Hence every indegree in~$\vec R$ is even. 
  Summing the indegrees in each component shows that every component of~$R$ has
  an even number of edges. 
  Since $\chi(R) \le 3$, \Cref{lem:partite-residual} gives $\lir(R) \le 10$.
  The edge-partition into~$P$, $O' \cup S$, and~$R$ now gives
  $\lir(G) \le 1 + 4 + 10 = 15$.

  For the second assertion, let~$G$ be decomposable. 
  Every component of~$G$ is then decomposable. 
  For each component having an odd number of edges, apply
  \Cref{lem:odd-reduction}, and let~$L$ be the union of the locally irregular
  subgraphs thereby obtained. 
  These subgraphs lie in distinct components of~$G$, so~$L$ is locally
  irregular and can be assigned at most one additional colour.
  Components of~$G$ having an even number of edges are unchanged, and the lemma
  ensures that every component left inside an originally odd component has an
  even number of edges. 
  Thus every component of $G-E(L)$ has an even number of edges. 
  Since $G-E(L)$ is still $6$-colourable, the first assertion gives at
  most~$15$ further colours.
  Hence $\lir(G) \le 1 + 15 = 16$.
\end{proof}

\section{Preserving connectivity}

\subsection{Selections}

For a graph~$H$, let $c(H)$ denote its number of components, including isolated
vertices.
For a connected graph~$H$ and $A \subseteq E(H)$, define
\begin{equation}\label{eq:capacity-formula}
  \mu_H(A) := |A| - c(H-A) + 1 \,.
\end{equation}
This is the maximum number of edges of~$A$ that can be deleted while preserving
connectivity.
Indeed, contract the components of~$H-A$. 
The resulting connected multigraph has $c(H-A)$ vertices, so at least
$c(H-A)-1$ edges of~$A$ must be retained. 
Conversely, retaining the edges corresponding to a spanning tree of this
multigraph and deleting all other edges of~$A$ attains this bound. 
The maximum-deletion interpretation also shows that~$\mu_H$ is monotone. 
In addition,
\begin{equation}\label{eq:capacity-uncrossing}
  \mu_H(A) + \mu_H(B)
    \ge \mu_H(A \cup B) + \mu_H(A \cap B)
    \qquad \text{for all } A,B \subseteq E(H) \, .
\end{equation}
Indeed, put $X := H-A$, $Y := H-B$, and $n := |V(H)|$.
Extend a spanning forest~$F_0$ of $X \cap Y$ to spanning forests~$F_X$ of~$X$
and~$F_Y$ of~$Y$.
The graph $F_X \cup F_Y$ has the same components as $X \cup Y$.
Since $F_0 \subseteq F_X \cap F_Y$,
\[
  n - c(X \cup Y) \le e(F_X \cup F_Y)
    \le e(F_X) + e(F_Y) - e(F_0)
    = n - c(X) - c(Y) + c(X \cap Y) \, .
\]
Thus $c(X) + c(Y) \le c(X \cap Y) + c(X\cup Y)$.
Since $X \cap Y = H - (A \cup B)$, $X \cup Y = H - (A \cap B)$, and $|A|+|B| =
|A \cup B| + |A \cap B|$, the inequality in \eqref{eq:capacity-uncrossing}
follows from~\eqref{eq:capacity-formula}.

\begin{lemma}\label{lem:connected-selection}
  Let~$H$ be connected, and let $(A_i)_{i \in I}$ be a finite family of subsets
  of~$E(H)$. 
  There are distinct edges $e_i \in A_i$, for $i \in I$, such that $H - \{e_i
  \colon i \in I\}$ is connected if and only if
  \begin{equation}
  \label{eq:selection-condition}
    \mu_H \left(\bigcup_{i \in J} A_i \right) \ge |J|
      \qquad \text{for every } J \subseteq I \, .
  \end{equation}
\end{lemma}
\begin{proof}
  Suppose first that distinct edges $e_i \in A_i$ can be chosen and that their
  simultaneous deletion leaves~$H$ connected. 
  For every~$J \subseteq I$, the $|J|$ edges in $\{e_i \colon i \in J\}$ belong
  to $\bigcup_{i \in J} A_i$ and may also be deleted while preserving
  connectivity. 
  The interpretation of~$\mu_H$ above gives~\eqref{eq:selection-condition}.

  Conversely, shrink the sets~$A_i$, one edge at a time, while
  preserving~\eqref{eq:selection-condition}. 
  Suppose that, in the resulting inclusion-minimal family, some~$A_k$ contains
  distinct edges $x,y$.
  Write $A(J) := \bigcup_{j \in J} A_j$.
  Removing~$x$ from~$A_k$ violates the condition, and the violating index set
  must contain~$k$. 
  Hence there is $J_x \subseteq I \setminus \{k\}$ such that, for
  $P := (A_k \setminus \{x\}) \cup A(J_x)$, we have $\mu_H(P) < |J_x| + 1$.
  Since~$\mu_H$ is integer-valued, $\mu_H(P)\le|J_x|$. 
  Similarly, there is $J_y \subseteq I \setminus \{k\}$ such that, for $Q :=
  (A_k \setminus \{y\}) \cup A(J_y)$, we have $\mu_H(Q) \le |J_y|$.
  Now $P \cup Q = A_k \cup A(J_x \cup J_y)$ and $P \cap Q \supseteq A(J_x \cap
  J_y)$.
  The selection condition and monotonicity give
  \[
    \mu_H(P \cup Q) \ge |J_x \cup J_y| + 1,
    \qquad
    \mu_H(P \cap Q) \ge |J_x \cap J_y| \, .
  \]
  Therefore, \eqref{eq:capacity-uncrossing} gives
  \[
    |J_x| + |J_y| \ge \mu_H(P) + \mu_H(Q)
                  \ge |J_x \cup J_y| + 1 + |J_x \cap J_y|
                  = |J_x| + |J_y| + 1 \, ,
  \]
  a contradiction.

  Thus every~$A_i$ is a singleton, say $A_i=\{e_i\}$. 
  Applying \eqref{eq:selection-condition} to~$I$ and using the interpretation
  of~$\mu_H$ gives
  \[
    |I| \le \mu_H \bigl(\{e_i \colon i \in I\}\bigr)
      \le |\{e_i \colon i \in I\}| \le |I| \, .
  \]
  Equality throughout shows that the edges~$e_i$ are distinct and that all of
  them can be deleted while preserving connectivity.
\end{proof}

\subsection{Selecting crossing edges while preserving connectivity}

\begin{lemma}\label{lem:switching}
  Let~$G$ be a connected graph, let $k\ge2$, and let~$\varphi$ be a
  crossing-optimal $k$-labelling of~$V(G)$. 
  Let~$C$ be its crossing graph.
  Then~$C$ is connected.
  Moreover, if $uv$ is monochromatic with $\varphi(u) = \varphi(v) = i$, then,
  for every label~$j\neq i$, the vertices $u$ and $v$ belong to the same
  component of $C[\varphi^{-1}(\{i,j\})]$.
\end{lemma}
\begin{proof}
  Let $uv$ be monochromatic with label~$i$, let $j\neq i$, and let~$W$ be the
  vertex set of the component of $C[\varphi^{-1}(\{i,j\})]$ containing~$u$.
  If $v\notin W$, interchange labels $i$ and~$j$ on~$W$. 
  The edge~$uv$ becomes crossing, and no crossing edge is lost: a crossing edge
  with exactly one end in~$W$ has its other end labelled outside~$\{i,j\}$.
  This contradicts crossing-optimality and proves the second assertion.
  In particular, the ends of every monochromatic edge belong to the same
  component of~$C$.
  Since~$G$ is connected, so is~$C$.
\end{proof}

For the next lemma and corollary, let~$G$ be a connected graph and let
$\varphi\colon V(G)\to[6]$ be crossing-optimal, with crossing graph~$C$ and
monochromatic graph~$F$. 
By \Cref{lem:switching}, $C$ is connected.
Fix a partition $B_1,B_2,B_3$ of~$[6]$ into pairs, write~$B(i)$ for the pair
containing~$i$, and let~$\vec F$ be any orientation of~$F$. Put
\[
  \mathcal I := \{(a,j) \colon a = (u,v) \in E(\vec F),\
                   j \in [6] \setminus B(\varphi(v))\} \, ,
\]
and, for each $(a,j) \in \mathcal I$, define
\[
  L(a,j) := \{vw\in E(C) \colon \varphi(w) = j\},
  \qquad \text{where } a = (u,v) \,.
\]

\begin{lemma}
\label{lem:support-capacity}
  For every $\mathcal X \subseteq \mathcal I$,
  \begin{equation}\label{eq:support-capacity}
    \mu_C \left(\bigcup_{(a,j) \in \mathcal X} L(a,j) \right) \ge |\mathcal X| \, .
  \end{equation}
\end{lemma}
\begin{proof}
  Fix~$\mathcal X \subseteq \mathcal I$, and put $A := \bigcup_{(a,j) \in
  \mathcal X} L(a,j)$.
  For each $(u,v) \in E(\vec F)$, \Cref{lem:switching} gives a path from~$u$
  to~$v$ in $C[\varphi^{-1}(B(\varphi(v)))]$. 
  Every edge of~$A$ has endpoint labels in two different sets
  among~$B_1,B_2,B_3$, whereas this path uses labels from only one of them.
  Thus the path lies in~$C-A$.

  Let~$\mathcal K$ be the set of components of~$C-A$. 
  Contract them, discard loops, and choose a spanning tree~$\Theta$ of the
  resulting connected multigraph, identifying its edges with their preimages
  in~$A$.
  Put $D := A \setminus E(\Theta)$.
  Then $C-D$ is connected and, by \eqref{eq:capacity-formula}, $|D| = \mu_C(A)$.

  Fix a pair~$p$ of labels lying in two distinct sets among~$B_1,B_2,B_3$, and
  identify its two labels with the elements of~$\mathbb Z_2$. 
  Let $\mathcal X_p$ consist of the pairs $(a,j) \in \mathcal X$ with
  $\{\varphi(v),j\} = p$, where $a = (u,v)$, and let~$D_p$ consist of the edges
  of~$D$ whose end labels form~$p$.
  Let~$\mathcal H$ be the set of heads occurring in~$\mathcal X_p$.
  Choose the variables $(\xi_v)_{v \in \mathcal H}$ independently and uniformly
  from~$\mathbb Z_2$, and put $\xi_v := 0$ for $v \notin \mathcal H$.

  Choose a root~$K_0 \in \mathcal K$ for~$\Theta$ and set $\sigma_{K_0} := 0$. 
  For an edge $KK'$ directed away from the root, represented by $xy$ with $x
  \in K$ and $y \in K'$, define recursively
  \begin{equation}\label{eq:tree-shifts}
    \sigma_{K'} - \sigma_K \equiv \xi_x - \xi_y \pmod 2 \, .
  \end{equation}
  These shifts are homogeneous linear forms in the variables~$\xi_v$.
  For every $K \in \mathcal K$ and every vertex~$v \in K$ whose label lies
  in~$p$, interchange the two labels of~$p$ precisely when $\xi_v + \sigma_K =
  1$ in~$\mathbb Z_2$.
  Leave all other labels unchanged.

  A crossing edge can be lost only if its original labels form~$p$. 
  Since $C-D = (C-A) \cup E(\Theta)$, first suppose that such an edge lies
  in~$C-A$.
  Neither end belongs to~$\mathcal H$, since otherwise the edge would lie
  in~$A$, and both ends receive the same component shift. 
  If instead it belongs to~$E(\Theta)$, the total shifts at its ends agree
  by~\eqref{eq:tree-shifts}. 
  Thus every edge of~$C-D$ remains crossing.

  For each $(a,j) \in \mathcal X_p$, write $a = (u,v)$. 
  The ends have the same original label and lie in the same component of~$C-A$,
  so their component shifts cancel. 
  Here $v \in \mathcal H$, and~$\xi_v$ is independent of~$\xi_u$ if $u \in
  \mathcal H$, while $\xi_u = 0$ otherwise. 
  Thus $\xi_v - \xi_u$ is uniform on~$\mathbb Z_2$, and the edge~$uv$ becomes
  crossing with probability~$1/2$.
  Each arc occurs at most once in~$\mathcal X_p$, so the expected gain is at
  least~$|\mathcal X_p|/2$.
  Only edges of~$D_p$ can be lost. 
  For each such edge, the final binary difference of its labels is~$1$ plus a
  homogeneous linear form in independent bits. 
  If the form is zero, this difference never vanishes; if it is nonzero, it is
  uniform on~$\mathbb Z_2$. 
  Thus the edge is lost with probability at most~$1/2$.
  By crossing-optimality and linearity of expectation,
  \[
    0 \ge \mathbb E[\Delta_{\mathrm{cr}}]
      \ge \frac{|\mathcal X_p| - |D_p|}{2} \, ,
  \]
  where $\Delta_{\mathrm{cr}}$ is the change in the number of crossing edges.
  As~$p$ ranges over all pairs of labels lying in two distinct sets
  among~$B_1,B_2,B_3$, the sets~$\mathcal X_p$ partition~$\mathcal X$ and the
  sets~$D_p$ partition~$D$.
  Summing $|\mathcal X_p| \le |D_p|$ gives $|\mathcal X| \le |D| = \mu_C(A)$.
\end{proof}

\begin{corollary}\label{cor:support}
  There are distinct edges $e_{a,j} \in L(a,j)$, for $(a,j) \in \mathcal I$,
  whose simultaneous deletion leaves~$C$ connected.
\end{corollary}
\begin{proof}
  Apply \Cref{lem:connected-selection} to the lists~$L(a,j)$, using
  \Cref{lem:support-capacity}.
\end{proof}

\subsection{Proof of the main theorem}

\begin{proof}[Proof of \Cref{thm:global18}]
  Since every component of a decomposable graph is decomposable and the same
  set of colours may be reused in distinct components, we may assume that~$G$
  is connected. 
  Choose a crossing-optimal $6$-labelling~$\varphi$ of~$V(G)$, and let~$C$
  and~$F$ be its crossing and monochromatic graphs. 
  If~$F$ is empty, then this labelling is proper, and the second assertion of
  \Cref{lem:even-fifteen} gives $\lir(G)\le16$.

  Otherwise, partition~$[6]$ into pairs~$B_1,B_2,B_3$, and put
  $V_i := \varphi^{-1}(B_i)$ and $J_i := F[V_i]$ for $i \in [3]$.
  The graphs $J_1,J_2,J_3$ edge-partition~$F$.
  For each~$J_i$ with an edge, choose an orientation~$\vec J_i$ as in
  \Cref{lem:general-core}, and take the empty orientation otherwise.
  Together these give an orientation~$\vec F$ of~$F$.
  Apply \Cref{cor:support}, and let~$P_i$ be the spanning subgraph of~$C$
  consisting of the selected edges~$e_{a,\ell}$ with $a \in E(\vec J_i)$ and
  $\ell \in [6] \setminus B_i$.
  These graphs are edge-disjoint, and each arc contributes one edge to each
  label outside~$B_i$.
  Thus~$P_i$ is bipartite between~$V_i$ and its complement, with
  \[
    d_{P_i}(x,V_j) = 2d^-_{\vec J_i}(x)
      \qquad \text{for every } x \in V_i,\ j \in [3] \setminus \{i\} \,.
  \]
  For every~$J_i$ with an edge, let~$Y_1,Y_2$ be the two sets in
  $\{V_j \colon j \in [3] \setminus \{i\}\}$ and apply \Cref{lem:general-core}
  to~$J_i$ with $X = V_i$ and $B = P_i$. 
  For either prescribed parity, this gives a subgraph~$Q_i \subseteq P_i$ for
  which $J_i \cup Q_i$ is locally irregular. 
  If~$J_i$ is edgeless, take~$Q_i$ empty.

  Since $F \neq \emptyset$, choose $i_0 \in [3]$ such that~$J_{i_0}$ has an edge.
  For each $i \neq i_0$ for which~$J_i$ has an edge, choose~$Q_i$ with an even
  number of edges. 
  Then choose~$Q_{i_0}$ with $e(Q_{i_0}) \equiv e(C) \pmod 2$. 
  Thus
  \[
    \sum_{i \in [3]} e(Q_i) \equiv e(C) \pmod 2 \,.
  \]
  Put $H_i := J_i \cup Q_i$ for every~$i \in [3]$.

  Put $R := C - \bigcup_{i \in [3]} E(Q_i)$. 
  The displayed congruence shows that~$R$ has an even number of edges and
  contains the connected spanning graph $C-\bigcup_{i \in [3]} E(P_i)$, which
  is connected by \Cref{cor:support}.
  Thus~$R$ is connected, and~$\varphi$ is a proper $6$-labelling of~$V(R)$.
  The graphs~$H_i$ with $E(J_i) \neq \emptyset$, together with~$R$,
  edge-decompose~$G$. 
  The first assertion of \Cref{lem:even-fifteen} gives $\lir(R) \le 15$, and
  therefore
  \[
    \lir(G) \le 3 + \lir(R) \le 3 + 15 = 18 \, .
  \]
\end{proof}

\section*{Acknowledgements}

This work was supported in part by the Coordena\c{c}\~{a}o de
Aperfei\c{c}oamento de Pessoal de N{\'i}vel Superior (CAPES), Finance Code 001,
and by the European Union's Horizon Europe Framework Programme under Grant
Agreement No.~101183743 (AGATE).
C.~N.~Lintzmayer was supported in part by the Conselho Nacional de
Desenvolvimento Cient{\'i}fico e Tecnol{\'o}gico (CNPq), grants 404315/2023-2
and 301143/2025-0, and by the Funda\c{c}\~{a}o de Amparo à Pesquisa do Estado
de S\~{a}o Paulo (FAPESP), grant 2024/16092-6.
G.~O.~Mota was supported by CNPq, grants 420838/2025-2 and 315916/2023-0, and
by FAPESP, grants 2024/13859-4 and 2023/03167-5.
M.~Sambinelli was supported in part by FAPESP, grant 2025/00323-1.
V.~F.~dos Santos was supported in part by CNPq, grants 308129/2025-3 and
404479/2023-5, and by FAPESP, grant 2025/27516-4.

\section*{Description of the use of artificial intelligence}

The development of this paper combined substantial use of ChatGPT (GPT-5.6 Sol)
with ideas and proofs contributed by the authors.
ChatGPT generated an initial draft presenting a 94-colour bound in response to
a simple prompt.
Subsequent improvements and simplifications resulted from an extended
mathematical interaction with the model and from the authors' mathematical
work.

After the initial long and complicated proof giving a 94-colour bound, we
provided mathematical guidance, proposed approaches, proved auxiliary lemmas,
and supplied material from several years of our research, including students'
theses and digital notes documenting unsuccessful attempts.
This interaction enabled the model to reduce the bound from 94 to 70. 
After polishing the arguments and improving some auxiliary lemmas, we reduced
the bound to 26 and then to 20. 
Further work with the model led to the 18-colour bound proved here.

The development therefore combined AI-generated suggestions with human
mathematical contributions. 
The research direction and the assessment of proposed arguments remained under
the authors' control.
The authors worked alongside ChatGPT, using it as an assistant to speed up
verification and provide ideas and intermediate results.

Beyond the initial 94-colour draft, the model assisted in exploring possible
arguments, organising the manuscript, and simplifying proofs and exposition. 
Its contribution was therefore not limited to language editing; however, the
final result should not be understood as an autonomous AI discovery.
The authors take full responsibility for the mathematical correctness,
attribution, and presentation of the results in this paper.

\bibliographystyle{amsplain}
\bibliography{main}

\end{document}